# Comparing two samples by penalized logistic regression

## Konstantinos Fokianos[*]

*Department of Mathematics & Statistics*
*University of Cyprus*
*e-mail:* fokianos@ucy.ac.cy

**Abstract:** Inference based on the penalized density ratio model is proposed and studied. The model under consideration is specified by assuming that the log–likelihood function of two unknown densities is of some parametric form. The model has been extended to cover multiple samples problems while its theoretical properties have been investigated using large sample theory. A main application of the density ratio model is testing whether two, or more, distributions are equal. We extend these results by arguing that the penalized maximum empirical likelihood estimator has less mean square error than that of the ordinary maximum likelihood estimator, especially for small samples. In fact, penalization resolves any existence problems of estimators and a modified Wald type test statistic can be employed for testing equality of the two distributions. A limited simulation study supports further the theory.

## 1. Introduction

The density ratio model is specified by assuming that the log–ratio of two unknown probability density functions is linear in some parameters, see Qin and Zhang (1997) and Qin (1998). The model is motivated by considering a logistic regression model for a binary random variable $Y$, which assumes the values 1 and 2–where "1" denotes success–and $X$ a $d$–dimensional vector of covariates, see Anderson (1972, 1979), Breslow and Day (1980), Cox and Snell (1989), for example. Then, the logistic regression model expresses the probability of the event $\{Y = 1\}$ as a function of $X$ by

$$P[Y = 1 \mid x] = \frac{\exp(\alpha^\star + \beta' x)}{1 + \exp(\alpha^\star + \beta' x)},$$

where $\alpha^\star$ is a scalar parameter and $\beta$ is a $d \times 1$ vector of regression coefficients. The logistic model leads to the density ratio model when considering

---

[*]Part of this work was carried out while the author was visiting the Department of Statistics, University of Munich. Discussions and support by L. Fahrmeir and K. Strimmer are greatly acknowledged. The comments of a reviewer improved the presentation.

case–control or retrospective sampling, see Prentice and Pyke (1979), Farewell (1979). Suppose that $X_{11}, \ldots, X_{1n_1}$ is a random sample from $G(x \mid Y = 1)$ and $X_{21}, \ldots, X_{2n_2}$ is another independent sample from $G(x \mid Y = 2)$. Set $\pi = \mathrm{P}[Y = 1]$ and $g(x \mid Y = i) = dG(x \mid Y = i)/dx$, for the conditional probability density function of $X$ given $Y = i$, $i = 1, 2$. Bayes' theorem shows that

$$
\begin{aligned}
\frac{g(x \mid Y = 1)}{g(x \mid Y = 2)} &= \frac{1 - \pi}{\pi} \exp\left(\alpha^\star + \beta' x\right) \\
&= \exp\left(\alpha + \beta' x\right),
\end{aligned}
$$

with $\alpha = \alpha^\star + \log\{(1 - \pi)/\pi\}$. The last equation justifies the term density ratio model: the densities of the observations are related by a parametric exponential tilt, but otherwise are unknown. In general, consider the following two independent samples semiparametric problem:

$$
\begin{aligned}
X_{11}, \ldots, X_{1n_1} &\text{ is a random sample from } & g_1(x) &= \exp\left(\alpha + \beta' h(x)\right) g_2(x), \\
X_{21}, \ldots, X_{2n_2} &\text{ is a random sample from } & g_2(x),
\end{aligned}
$$

$$(1)$$

where $g_i(x)$, $i = 1, 2$ are unknown probability density functions. The quantity $\alpha$ is an unknown scalar, $\beta$ is a $d$–dimensional vector of parameters while $h(x)$ is a $d$-dimensional vector which consists of known functions of $X$. Several examples of distribution fall within the above framework, in particular the exponential family of distributions satisfies (1) trivially. An important observation is that when the model holds, and if $\beta = 0$, then the two samples are identically distributed. We conclude that model (1) is useful to the semiparametric comparison of two samples in the sense that the densities $g_i(.)$, $i = 1, 2$ are left completely unspecified but the weight function $\exp\left(\alpha + \beta' h(x)\right)$ depends on some finite dimensional parameter. The last remark connects the density ratio model and biased sampling theory, see Vardi (1982, 1985), Gill et al. (1988), Gilbert et al. (1999) and Gilbert (2000).

Inference regarding both finite and infinite dimensional parameters of model (1) has been studied by various authors assuming that the sample sizes tend to infinity in a suitable way. The methodology is based on empirical likelihood inference, see Owen (2001). Accordingly, a parametric likelihood function for the finite dimensional parameters is obtained after profiling out the infinite dimensional parameter of the model. However, there are applications where the sample sizes are small and therefore direct application of the empirical likelihood methodology might not be suitable. To overcome these problems we put forward a penalized empirical likelihood function for inference, see (4), which depends on two forms, a regularization parameter and a penalty function. The regularization parameter controls the amount of penalization whereas the penalty term is a function of the finite dimensional parameter of the model. The choice of the penalty function is independent of the infinite dimensional parameter as it will be explained in Section 2 where the methodology is developed in detail. Using penalized empirical likelihood leads to solution of several problems including elimination of non existence and non convergence issues. In addition approximations even for relatively even small sample sizes are adequate enough to allow for testing. Section 3 reviews several facts about the resulting estimators and shows that a judicious choice of the regularization parameter leads to consistent and asymptotically normal estimators. The last section examines the finite sample performance of the proposed estimators by some limited simulations. The paper closes with several remarks and an appendix.

## 2. Penalized likelihood inference

Suppose that $X_{ij}$, $j = 1, \ldots, n_i$, $i = 1, 2$ are two independent random samples from cumulative distribution functions $G_i(x)$, $i = 1, 2$, respectively. Let $g_i(x)$ be the corresponding density functions, that is $g_i(x) = dG_i(x)/dx$, and denote by $\mathbf{X} = (X_{11}, X_{12}, \ldots, X_{1n_1}, X_{21}, \ldots, X_{2n_2})'$ the vector of all observations. Suppose further that the density ratio model (1) holds. Inference for the finite dimensional parameters $\theta = (\alpha, \beta')'$ is based on the so–called empirical likelihood–see Owen (2001). Accordingly, let $p_{ij}$ denote the size of the jump at the observed datum $x_{ij}$, that is $p_{ij} = dG_2(x_{ij}) = G_2(x_{ij}^+) - G_2(x_{ij}^-)$, $j = 1, 2, \ldots, n_i$, $i = 1, 2$ and consider the following nonparametric likelihood given the data,

$$
\begin{aligned}
L(\alpha, \beta, G_2 \mid \mathbf{x}) &= \left\{ \prod_{j=1}^{n_1} \exp\left(\alpha + \beta' h(x_{1j})\right) dG_2(x_{1j}) \right\} \left\{ \prod_{j=1}^{n_2} dG_2(x_{2j}) \right\} \\
&= \left\{ \prod_{i=1}^{2} \prod_{j=1}^{n_i} p_{ij} \right\} \left\{ \prod_{j=1}^{n_1} \exp\left(\alpha + \beta' h(x_{1j})\right) \right\}.
\end{aligned} \tag{2}
$$

Following Qin and Zhang (1997) and Fokianos et al. (2001), elimination of the infinite dimensional parameter $G_2(.)$ is accomplished by maximizing the first product of (2) subject to the constraints $\sum_{i=1}^{2} \sum_{j=1}^{n_i} p_{ij} = 1$, and $\sum_{i=1}^{2} \sum_{j=1}^{n_i} p_{ij}(\exp(\alpha + \beta' h(x_{ij}) - 1)) = 0$. Avoiding unnecessary repetition, the empirical log-likelihood is given by

$$
l(\theta) = -\sum_{i=1}^{2} \sum_{j=1}^{n_i} \log\left[1 + \rho_1 \exp(\alpha + \beta' h(x_{ij}))\right] + \sum_{j=1}^{n_1} (\alpha + \beta' h(x_{1j})) \tag{3}
$$

where $\rho_1 = n_1/n_2$. Furthermore, if $\hat{\theta} = (\hat{\alpha}, \hat{\beta})'$ denotes the maximum likelihood estimator of $\theta$, assuming that it exists, then it can be shown that

$$
\hat{p}_{ij} = \frac{1}{n_2} \frac{1}{1 + \rho_1 \exp(\hat{\alpha} + \hat{\beta}' h(x_{ij}))}.
$$

Hence a consistent estimator for both of $G_1(.)$ and $G_2(.)$ can be constructed provided that the total sample size $n = n_1 + n_2$ tends to infinity such that $n_1/n_2 \to \rho_1$–see Qin and Zhang (1997) and Fokianos et al. (2001) for further details.

Expression (3), after reparametrization, is equivalent to standard logistic regression likelihood–a direct consequence of the equivalence between retrospective and prospective sampling, see Prentice and Pyke (1979). Accordingly, inference for the parameter $\theta$ is carried out by numerous statistical programs which include logistic regression modeling. However, it is well known that if the sample size is large, then the logistic likelihood equations yield the maximum likelihood estimator provided that it exists. In this case, the score equations are solved by the Fisher scoring method and under some reasonable regularity assumptions a sequence of approximations is derived that converges to the maximum likelihood estimators of $\alpha$ and $\beta$. In contrast, convergence and existence issues arise when the sample size is relatively small. When the sample size is small, maximum likelihood estimates fail to exist in general if the data are completely, quasi–completely or partially separated. More specifically, there does not exist a maximum likelihood estimator of $\beta$ if there exists a hyperplane in the covariate space such that when $Y = 1$, the covariates belong to the left side of the hyperplane whereas for $Y = 2$, the covariates belong to the other side of the hyperplane (see Albert and Anderson (1984), McCullagh and Nelder (1989, Sec. 4.4) and Santner and Duffy (1989, p. 234)). The discussion shows that application of the density ratio model is questionable when only a few observations belong to each group.

However, regularization of the log –likelihood function (3), in the sense of adding a concave penalty term, avoids existence and convergence issues, see also Hastie and Tibshirani (2004). Consider the so–called penalized empirical log-likelihood function

$$l_p(\theta) = l(\theta) - \lambda_n J(\beta). \tag{4}$$

where $l(\theta)$ denotes the unrestricted log–likelihood function given by (3), $\lambda_n$ is a sequence of regularization parameter controlling the amount of shrinkage and $J(.)$ is a penalty function on the parameter $\beta$. The parameter $\alpha$–the intercept– is not penalized explicitly because it is a function of $\beta$ when (1) holds. It is well known that when $\lambda_n \to 0$ then (4) yields to the unrestricted maximum likelihood estimator whereas if $\lambda_n \to \infty$ then $\theta$ shrinks towards 0.

We argue that it is reasonable to employ the profile empirical likelihood for inference. Expression (4) is not a proper log–likelihood function in the sense that the first summand, namely the quantity $l(\theta)$, is the outcome of a profiling procedure derived by means of empirical likelihood methodology. However, recent work in the area of semiparametric statistical inference shows that such functions share the same properties of common likelihood functions, see Murphy and van der Vaart (2000). Thus, it is sensible to penalize (3) by introducing an extra concave term which does not depend on the distribution function $G_2$. The method resolves successfully the existence problem of the maximum likelihood estimator of $\theta$ and yields estimates of the unknown distribution functions even for relatively small samples.

Motivation of penalized empirical likelihood (4) stems also from a Bayesian point of view. Suppose that the density ratio model (1) holds conditionally on

$\beta$ and suppose that $\beta$ is a random variable with a density proportional to

$$\pi(\beta) \propto \exp\left(-\lambda_n J(\beta)\right).$$

The above display points out to a crucial point, namely the independence between the prior of $\beta$ and the cdf $G_2$. It is precisely this fact that enables to base inference on the posterior likelihood function

$$\pi(\beta|G_2, \mathbf{x}) \propto L(\alpha, \beta, G_2|\mathbf{x})\pi(\beta).$$

By setting $l_p(\theta) \equiv \log \pi(\beta|G_2, \mathbf{x})$ we obtain the penalized log–likelihood expression (4). Indeed, maximization of the last display with respect to $p_{ij}$'s is equivalent to maximization of (2) since the second factor does not depend on the infinite dimensional parameter. That is *the penalized log–likelihood function still produces valid cdf estimators that satisfy the constraints imposed by the density ratio model*. We conclude that maximization of $\pi(\beta|G_2, \mathbf{x})$ has no effect on the final form of the profile log–likelihood function as long as the designated function of $\beta$ that multiplies (2) does not depend on $G_2$. Note that the methodology is quite general and it is rather interesting to study its properties for a wider class of examples associated with empirical likelihood theory.

The penalized log–likelihood equations (4) depends on the choice of regularization parameter and the penalty term $J(\beta)$. The selection of $\lambda_n$ can be based on existing cross–validation methods but it is unclear what are the properties of such an approach especially when few observations are available. There are various approaches to choose the penalty function–perhaps the most popular being $J(\beta) = \sum_{j=1}^{d} \beta_j^2$ leading to the so called ridge regression type estimators (see Le Cessie and Van Houwelingen (1992) for the case of the logistic regression when the number of covariates is large). A more general family of penalty functions is given by $J(\beta) = \sum_{j=1}^{d} \gamma_k \psi(\beta_j)$, with $\gamma_k > 0$ and leads to several well known examples. For instance, if $\gamma_k = 1$ and $\psi(\beta_j) = |\beta_j|$, then $J(\beta)$ reduces to the $L^1$ penalty (see Tibshirani (1996)) while the $L^q$ penalty is obtained when $\psi(\beta_j) = |\beta_j|^q$, for $q > 0$, see Frank and Friedman (1993). Insight on the choice of penalty function is given in the recent articles of Antoniadis and Fan (2001) and Fan and Li (2001). However, in what follows consider the penalty function

$$J(\beta) = \sum_{j=1}^{d} |\beta_j|^q, \quad q > 1. \tag{5}$$

In particular, the ridge penalty will be employed in applications since the linear combination $\beta' h(x)$ is expected to vary smoothly in accordance with (1).

## 3. Main results

To fix notation, suppose that $\hat{\theta} = (\hat{\alpha}, \hat{\beta}')'$ denotes the maximum likelihood estimator derived by maximization of the unrestricted log–likelihood (3) and let $\hat{\theta}^\lambda$ denote the constrained maximum likelihood estimator obtained by maximizing

(4). In what follows, $\theta_0$ denotes the true value of the parameter $\theta$ when (1) holds. Define

$$A_{11}(t) = \int_{-\infty}^{t} \frac{\exp(\alpha + \beta'h(x))}{1 + \rho_1 \exp(\alpha + \beta'h(x))} dG_2(x), \qquad A_{11} = A_{11}(\infty),$$

$$A_{21}(t) = \int_{-\infty}^{t} \frac{\exp(\alpha + \beta'h(x))}{1 + \rho_1 \exp(\alpha + \beta'h(x))} h(x) dG_2(x), \qquad A_{21} = A_{21}(\infty),$$

$$A_{22}(t) = \int_{-\infty}^{t} \frac{\exp(\alpha + \beta'h(x))}{1 + \rho_1 \exp(\alpha + \beta'h(x))} h(x) h'(x) dG_2(x), \quad A_{22} = A_{22}(\infty),$$

and put

$$A = \left( \begin{array}{cc} A_{11} & A_{21} \\ A_{21}' & A_{22} \end{array} \right), \quad V = \frac{\rho_1}{1 + \rho_1} A - \rho_1 \left( \begin{array}{c} A_{11} \\ A_{12}' \end{array} \right) \left( \begin{array}{cc} A_{11}, & A_{12} \end{array} \right).$$

The following lemma summarizes some useful facts regarding the large sample behavior of both the score function and the Hessian matrix. It is used for the proofs of Theorems 3.1, 3.2 and 3.3 and it is stated for completeness of the presentation.

**Lemma 3.1** (Qin and Zhang (1997)) *Suppose that the density ratio model (1) holds and assume that A is nonsingular. Denote by*

$$\nabla l(\theta) = \left( \frac{\partial l(\theta)}{\partial \alpha}, \frac{\partial l(\theta)}{\partial \beta} \right)',$$

*a $(d+1)$-dimensional vector. Then*

- *The score function of the unrestricted likelihood (3) is asymptotically normal,*

$$n^{-1/2} \nabla l(\theta_0) \to \mathcal{N}_{d+1}(0, V),$$ (6)

  *as $n \to \infty$ such that $n_1/n_2 \to \rho_1$.*

- *The Hessian matrix of the unrestricted likelihood (3) converges in probability, as $n \to \infty$ such that $n_1/n_2 \to \rho_1$, to*

$$-\frac{1}{n} \nabla^2 l(\theta_0) \to \frac{\rho_1}{1 + \rho_1} A \equiv S.$$ (7)

The large sample behavior of the score and the Hessian matrix requires a number of standard regularity conditions to be satisfied. However, in the context of logistic regression, these conditions can be easily verified, see for instance Santner and Duffy (1989). Hence, the density ratio model can be applied to a variety of settings provided that the log likelihood ratio of two probability densities whose support is identical, is a sufficiently smooth function. Several well known examples of distributions fall within this class and therefore model (1) is applicable to a large collection of problems.

The first task is to establish existence of $\hat{\theta}^\lambda$, that is whether the problem of maximizing (4) is well defined and if the resulting estimator is unique. When considering penalty functions of the form (5)–and more generally sums of smooth convex functions–the answer is given by the following result:

**Lemma 3.2** (Fu (1998)) *For the penalty function (5) and for given $q > 1$ and $\lambda_n > 0$ there exists a unique solution of the penalized score equations*

$$\nabla l_p(\theta) = \nabla l(\theta) - \lambda_n \nabla J(\beta).$$

*The unique solution is equal to the unique estimator of the maximization problem $\max_\theta l_p(\theta)$, provided that $\nabla l(\theta)$ is continuously differentiable with respect to $\theta$ and $\nabla^2 l(\theta)$ is positive semi definite.*

The equivalence of the density ratio model to logistic regression implies that both conditions of Lemma 3.2 are satisfied, or, in other words, a unique restricted estimator $\hat{\theta}^\lambda$ exists. Summarizing, the issue of non existence of maximum likelihood estimators in the context of logistic regression is successfully resolved by penalization. Thus, the density ratio model is applicable to situations where the available sample size is small. Furthermore, it will be shown that the choice of $\lambda_n = O(\sqrt{n})$ yields to $\sqrt{n}$–consistent estimators. In the case of $\lambda_n = o(\sqrt{n})$ then $\hat{\theta}^\lambda$ is consistent. These conclusions are supported further by some limited simulation results which are reported in the next section.

**Theorem 3.1** *Suppose that the density ratio model (1) holds and assume the regularity conditions of Lemma 3.1. Assume that the true parameter vector $\theta_0$ lies in a compact set. Suppose further that the penalty function is given by (5) and $\lambda_n/\sqrt{n} \to \lambda_0 \geq 0$. If $n \to \infty$ in such a way that $n_1/n_2 \to \rho_1$, then there exists a unique maximizer $\hat{\theta}^\lambda$ of (4) such that $\| \hat{\theta}^\lambda - \theta_0 \| = O_p(n^{-1/2})$. In particular, if $\lambda_0 = 0$, then $\hat{\theta}^\lambda$ is consistent.*

In addition, estimates of the unknown distribution functions are computed by setting

$$\hat{p}_{ij}^\lambda = \frac{1}{n_2} \frac{1}{1 + \rho_1 \exp(\hat{\alpha}^\lambda + \hat{\beta}^{\lambda'} h(x_{ij}))}.$$

Then

$$\hat{G}_2^\lambda(x) = \sum_{i=1}^2 \sum_{j=1}^{n_i} \hat{p}_{ij}^\lambda I(X_{ij} \leq x), \qquad (8)$$

and

$$\hat{G}_1^\lambda(x) = \sum_{i=1}^2 \sum_{j=1}^{n_i} \exp(\hat{\alpha}^\lambda + \hat{\beta}^{\lambda'} h(x_{ij})) \hat{p}_{ij}^\lambda I(X_{ij} \leq x),$$

where $I(.)$ is the indicator function. The jump sizes $\hat{p}_{ij}^\lambda$ sum up to 1 by construction since the penalized log likelihood function (4) is employed after profiling out the parameters $p_{ij}$.

**Remark 3.1** Identifiability of $\alpha$, $\beta$ and $G_2(.)$ is guaranteed by the work of Gilbert et al. (1999, Th. 2) which shows that if there exist a value $x_0$ such that $h(x_0) = 0$, then the density ratio model (1) is identifiable.

Next, consider the asymptotic distribution of $\hat{\theta}^\lambda$. It is useful to introduce the quantities

$$b = \nabla J(\beta), \tag{9}$$

a $d + 1$ dimensional vector, and

$$D = \nabla^2 J(\beta) \tag{10}$$

a $(d + 1) \times (d + 1)$ diagonal matrix. The choice of (5) yields

$$b = (0, \text{sign}(\beta_1)q \mid \beta_1 \mid^{q-1}, \dots, \text{sign}(\beta_d)q \mid \beta_d \mid^{q-1})',$$

and

$$D = \text{diag}(0, q(q-1) \mid \beta_1 \mid^{q-2}, \dots, q(q-1) \mid \beta_d \mid^{q-2})'.$$

**Theorem 3.2** *Suppose that the conditions of Theorem 3.1 are satisfied and let $\lambda_n$ be such that $\lambda_n/\sqrt{n} \to \lambda_0 \geq 0$. Then the unique maximizer $\hat{\theta}^\lambda$ tends to a $(d + 1)$–dimensional normal distribution. That is*

$$\sqrt{n}\big(\hat{\theta}^\lambda - \theta_0 + \lambda_0 S^{-1} b\big) \to \mathcal{N}_{d+1}(0, \Sigma),$$

*in distribution, as $n \to \infty$, where the $(d + 1) \times (d + 1)$ matrix $\Sigma = S^{-1} V S^{-1}$ is defined by means of (6) and (7), and $b$ is given by (9). In particular, if $\lambda_0 = 0$, then*

$$\sqrt{n}(\hat{\theta}^\lambda - \theta_0) \to \mathcal{N}_{d+1}(0, \Sigma).$$

The above theorem, when compared to the results of Lemma 3.1, does not point out to any advantages of the penalization for the density ratio model. Indeed, under the aforementioned conditions the maximum penalized empirical likelihood estimator of $\theta$ is asymptotically biased and its covariance matrix is equal to the covariance matrix of the unrestricted maximum empirical likelihood estimator. However the results is asymptotic and a careful examination of the proof of the theorem 3.2 shows that the asymptotic covariance matrix $\Sigma$ is approximated by

$$\hat{\Sigma} = \left(\frac{\rho_1}{1+\rho_1}\hat{A} + \frac{\lambda_n}{n}\hat{D}\right)^{-1} \hat{V} \left(\frac{\rho_1}{1+\rho_1}\hat{A} + \frac{\lambda_n}{n}\hat{D}\right)^{-1}, \tag{11}$$

where $\hat{A}$ is the empirical estimator of $A$ and $\hat{D}$ is equal to the diagonal matrix $D$ evaluated at $\hat{\beta}^\lambda$, provided that the penalty function is twice differentiable. Hence, for a judicial chosen value of $\lambda_n$, the mean square error of $\hat{\theta}^\lambda$ will be smaller than that of $\hat{\theta}$ in small samples. This point is further illustrated in the next section by considering some finite sample properties of the restricted maximum likelihood estimator, see equation (13) and subsequent discussion.

The estimate (11) performs satisfactorily even for small sample sizes given a known value of the regularization parameter–see next section. For large $n$, and if $\lambda_n = o(\sqrt{n})$, formula (11) reduces to that used by Qin and Zhang (1997) and Fokianos et al. (2001) for the asymptotic variance estimator of the unrestricted maximum likelihood $\hat{\theta}$.

Theorem 3.2 suggests that a standard Wald test can be used to test the hypothesis $H_0 : \beta = 0$–that is the two samples are identically distributed. In particular, under the hypothesis and provided that $\lambda_n = o(\sqrt{n})$, we obtain that the test statistic

$$\mathcal{W} = n\hat{\beta}^{\lambda'}\hat{\Sigma}_{22}^{-1}\hat{\beta}^{\lambda} \tag{12}$$

where $\hat{\Sigma}_{22}^{-1}$ denotes the estimated asymptotic variance of $\hat{\beta}^{\lambda}$ which is obtained by means of (11). Under the null hypothesis, the asymptotic distribution of $\mathcal{W}$ is the chi–square with $d$ degrees of freedom provided that the regularization parameter satisfies the aforementioned conditions. It is clear that (12) depends on the choice of the regularization parameter but the notation is suppressed for ease of presentation. A limited simulation study–see Section 4– shows that the chi–square approximation is valid even though the sample size is relatively small.

**Remark 3.2** Application of the density ratio model relies on the assumed functional form of $h(.)$. Misspecification of this form results in biased estimators and loss of efficiency– the point is made by Fokianos and Kaimi (2006) who employ the Box–Cox family of transformations to estimate $h(.)$. However, we have tacitly assumed throughout the presentation that the function $h(.)$ is the true function associated with the density ratio model (1).

**Remark 3.3** An alternative penalization scheme is based on (5) but with $0 < q \le 1$, see Knight and Fu (2000) for example, who studied the problem in the linear regression setup. These penalty functions are rather appealing in the sense that they combine model selection and estimation. Theory regarding the density ratio model in connection to penalty functions of the form $J(\beta)$ when $0 < q \le 1$ might be particular useful especially for multivariate observations. In this case, the main task is to identify a parsimonious functional form of the log likelihood ratio of two, or more, multivariate probability density functions and estimate some of the corresponding coefficients while setting the remaining to zero. In other words, the problem of modeling multivariate data by the density ratio model in large dimensions reduces to that of a model selection and estimation problem.

The last part of this section is devoted to the study of the large sample properties of the estimated cdf $G_2$ given by (8). The following theorem studies its asymptotic distributions and generalizes the corresponding theorem of Qin and Zhang (1997).

**Theorem 3.3** *Suppose that the conditions of Theorem 3.1 are satisfied and let $\lambda_n$ be such that $\lambda_n/\sqrt{n} \to \lambda_0 \ge 0$. Then, as $n \to \infty$ in such a way that $n_1/n_2 \to \rho_1$,*

$$\sqrt{n}\left(\hat{G}_2^{\lambda}(t) - \hat{G}_2(t) + \rho_1\lambda_0(A_{11}(t), A_{12}(t))S^{-1}b\right) \to W, \qquad in \quad D[-\infty, \infty]$$

*weakly, where $W$ is a Gaussian process possessing continuous sample paths and has mean 0 and covariance function specified by equation (15).*

A similar large sample result holds for $\hat{G}_1^\lambda$. Its proof is along the lines of the previous theorem and is omitted. In the next section we provide some empirical evidence for the finite sample performance of the estimate $\hat{\theta}^\lambda$.

## 4. Simulation study

The above results show that the penalized density ratio model depends on both the choice of the regularization parameter and the selection of the penalty function $J(\beta)$. In the sequel we assume that $J(\beta) = \sum_{i=1}^{d} \beta_i^2$ relying on the fact that the linear combination $\beta' h(x)$ is a smooth function $\beta$. This choice leads to the so–called ridge regression type estimators whose mean square error is less than that of the corresponding unrestricted maximum likelihood estimators, for a range of values of $\lambda_n$. For the ridge penalty, (see, for instance, Le Cessie and Van Houwelingen (1992)) consider the following approximate expansion of the penalized likelihood function for

$$\nabla l_p(\hat{\theta}^\lambda) \approx \nabla l_p(\theta_0) + (\hat{\theta}^\lambda - \theta_0)' \nabla^2 l_p(\theta_0).$$

After rearranging terms and taking into account equations (9), (10) and Lemma 3.1, an approximate expression for $\hat{\theta}^\lambda$ is given by

$$\hat{\theta}^\lambda \approx \left(\nabla^2 l(\theta_0) + 2\lambda_n D_r\right)^{-1} \left(\nabla^2 l(\theta_0)\theta_0 - \nabla l(\theta_0)\right),$$

where

$$D_r = \begin{bmatrix} 0 & 0 \\ 0 & I_d \end{bmatrix},$$

where $I_d$ denoted the $d$–dimensional identity matrix. On the other hand, a similar argument shows that $\theta_0 = \hat{\theta} + \{\nabla^2 l(\theta_0)\}^{-1}\nabla l(\theta_0)$. These two displays when combined show that

$$\hat{\theta}^\lambda \approx \left(\nabla^2 l(\theta_0) + 2\lambda_n D\right)^{-1} \nabla^2 l(\theta_0)\hat{\theta}. \tag{13}$$

Therefore $\hat{\theta}^\lambda$ shrinks towards zero as the value of the regularization parameter increases while for $\lambda_n \to 0$, $\hat{\theta}^\lambda$ is approaching the unrestricted maximum likelihood estimator. In particular, the above representation forms the basis for developing an asymptotic expression of the mean square error of $\hat{\theta}^\lambda$. The methodology is quite analogous to ordinary ridge regression and therefore is omitted (for more see Hoerl and Kennard (1970b,a). Consequently, it can be shown that the mean square error of $\hat{\beta}^\lambda$ attains its minimum value at some value of the regularization parameter, see the lower left hand side plot of Figure 1 for an example. In other words, the maximum penalized empirical likelihood estimator attains smaller mean square error than that of the estimator proposed by Qin and Zhang (1997), especially for small samples.

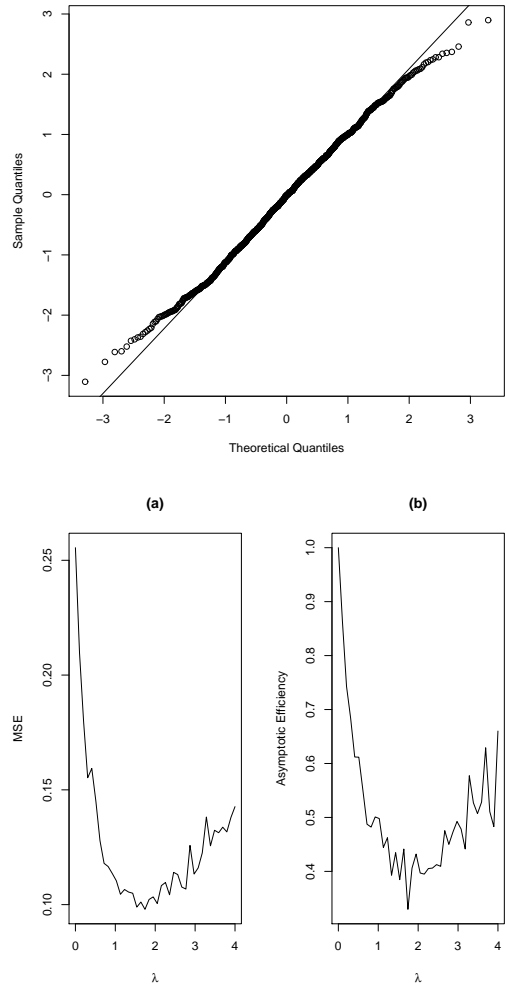

FIG 1. *Top: QQ–plot of the test statistic (12) under the hypothesis $\beta = 0$ for $\lambda = 1$. Here the data have generated according two lognormal populations with $\mu_1 = \mu_2 = 0$, $\sigma_1 = \sigma_2 = 1$ and $n_1 = n_2 = 10$. Bottom: (a) Mean square error of $\hat{\beta}^\lambda$ and (b) asymptotic efficiency of $\hat{\beta}$ with respect to $\hat{\beta}^\lambda$ as functions of $\lambda$. The data are drawn from two lognormal populations with $\mu_1 = 1$, $\mu_2 = 0$ $\sigma_1 = \sigma_2 = 1$, and $n_1 = n_2 = 20$. The quadratic penalty function has been used for fitting the model and the results are based on 1000 simulations.*

| True $\beta$ | Sample Size | $\lambda_n$ | $\beta^\lambda$ | MSE($\beta^\lambda$) | Power |
|---|---|---|---|---|---|
| 0 | $n_1 = n_2 = 10$ | 0 | 0.0193 | 0.3488 | 0.021 |
| | | 0.5 | 0.0089 | 0.2301 | 0.043 |
| | | 1.0 | -0.0181 | 0.1821 | 0.050 |
| | $n_1 = 10$ | 0 | -0.0189 | 0.1728 | 0.038 |
| | $n_2 = 30$ | 0.5 | 0.0070 | 0.1540 | 0.057 |
| | | 1.0 | 0.0027 | 0.1138 | 0.045 |
| 1 | $n_1 = n_2 = 10$ | 0 | 1.4123 | 9.7126 | 0.402 |
| | | 0.5 | 1.0023 | 0.2486 | 0.530 |
| | | 1.0 | 0.8525 | 0.1718 | 0.550 |
| | $n_1 = 10$ | 0 | 1.1321 | 0.3999 | 0.707 |
| | $n_2 = 30$ | 0.5 | 1.0122 | 0.1746 | 0.750 |
| | | 1.0 | 0.8978 | 0.1450 | 0.749 |

These results are empirically verified by considering the following simple example. Consider the case of two lognormal random samples with corresponding density functions

$$g(x, \mu_i, \sigma^2) = \frac{1}{x\sigma\sqrt{2\pi}} \exp\left(-(\log x - \mu_i)^2/2\sigma^2\right), \quad x > 0$$

for $i = 1, 2$. The density ratio model (1) is satisfied with $\alpha = (\mu_2^2 - \mu_1^2)/2\sigma^2$, $\beta = (\mu_1 - \mu_2)/\sigma^2$ and $h(x) = \log x$ and consider penalized empirical likelihood inference based on (4) by choosing appropriately both the penalty parameter and the penalty function. In this limited simulation study, the parameter $\lambda_n$ was chosen as a constant satisfying therefore the conditions of Theorems 3.1 and 3.2. Furthermore, a simple quadratic penalty function, that is $J(\beta) = \beta^2$, was used for maximizing the penalized likelihood (4). Table 1 summarizes the results of 1000 runs for different sample sizes and for values of $\lambda_n$. The first six lines of Table 1 report results when both groups of data are generated by the same log–normal distribution with $\mu_1 = \mu_2 = 0$ and $\sigma_1 = \sigma_2 = 1$. The mean square error of the penalized estimator of $\beta$, say $\hat{\beta}^\lambda$ is less than the mean square error of the unrestricted estimator even for small sample sizes. In addition, the nominal significance level of 5% for testing $\beta = 0$ by using the test statistic (12) is achieved satisfactorily in all cases except that of $\lambda_n = 0$. The variance estimate is obtained by means of (11). The last six lines of Table 1 report results under the alternative hypothesis. In this case the two log–normal samples were generated with $\mu_1 = 1$, $\mu_2 = 0$ and $\sigma_1 = \sigma_2 = 1$. The results are along the lines of the previous findings. It is interesting to observe that for $n_1 = n_2 = 10$, the mean square error of the unrestricted maximum likelihood estimator is large when compared with the mean square error of the penalized maximum likelihood estimators.

The top plot of Figure 1 shows a QQ–plot of the test static (12) under the hypothesis $H_0 : \beta = 0$. Notice that the approximation is quite satisfactory even though the sample sizes are relatively small. The bottom left hand side plot of Figure 1 shows the mean square error of $\hat{\beta}^\lambda$ as a function of $\lambda$. As the plot illustrates, there exists a value of $\lambda$ such that the mean square error is minimized–recall the above discussion. The right hand side of the plot shows the asymptotic efficiency of $\hat{\beta}$ with respect to $\hat{\beta}^\lambda$ defined by

$$e(\hat{\beta}; \hat{\beta}^\lambda) = \frac{\text{MSE}[\hat{\beta}^\lambda]}{\text{MSE}[\hat{\beta}]}.$$

For this specific example, the graph illustrates empirically that $e(\hat{\beta}; \hat{\beta}^\lambda) \leq 1$, implying that the restricted estimator is more efficient than the unrestricted estimator.

## 5. Conclusions

The density ratio model is a semiparametric alternative to the problem of comparing two, or more distributions functions. However, there are several drawbacks for its application especially when small samples are under consideration. The first part of this work dealt explicitly with the penalization approach to the empirical likelihood methodology. It was shown that the method is also motivated by Bayesian arguments and it is worth considering its performance to other settings. As a result, a new two sample methodology emerges where estimation of the parametric part can be carried out by a number of statistical programs while estimation of the non–parametric part is achieved by simple calculations. In particular the problem of existence of maximum empirical likelihood estimators is resolved successfully. In addition, the final estimators have smaller mean square error than the estimators without penalty. The test statistic–see (12)–for testing the equality of two distributions is obtained and the results were applied to simulated case control data. Application of the density ratio model depends on a number of assumptions, in particular the assumed functional form of model (1), the regularization parameter and the form of the penalty function (5). It is suggested to vary the regularization parameter within a predetermined range of values and use different functions $h(x)$ to examine the sensitivity of the results. Penalization by the quadratic penalty seems to be a sensible idea for several applications. As a final remark, we note that the problem can be generalized further by considering multiple biased sampling models–or semiparametric ANOVA–as in Fokianos et al. (2001). This will allow for semiparametric comparisons among several treatments and will also yield estimates of the unknown cdf. Furthermore, it is worth studying estimation of functionals associated with the baseline distribution, see Zhang (2000), in the presence of a regularization parameter. Different penalty functions can be taken into account and this methodology is quite promising for semiparametric comparison of multivariate distribution where the true functional form of (1) is complicated.

## Appendix

### Proof of Theorem 3.1

Let $u$ be a vector with $\| u \| = M$ where $M$ is a large constant. Suppose that for every $\epsilon > 0$ there exists $M > 0$ such that

$$\mathrm{P}\Big\{ \sup_{\|u\|=M} l_p(\theta_0 + n^{-1/2}u) < l_p(\theta_0) \Big\} \geq 1 - \epsilon. \tag{14}$$

That means there exists a local maximizer with probability tending to 1, denoted by $\hat{\theta}^\lambda$, in the ball $\{\theta_0 + n^{-1/2}u : \| u \| \leq M\}$ such that $\| \hat{\theta}^\lambda - \theta_0 \| = O_p(n^{-1/2})$. Now

$$
\begin{aligned}
l_p(\theta_0 + n^{-1/2}u) - l_p(\theta_0) &= \Big( l(\theta_0 + n^{-1/2}u) - l(\theta_0) \Big) \\
&\quad - \lambda_n \Big( J(\beta_0 + n^{-1/2}u) - J(\beta_0) \Big) \\
&= E_1 + E_2
\end{aligned}
$$

But

$$
E_1 = n^{-1/2}\nabla' l(\theta_0)u + \frac{1}{2}n^{-1}u'\nabla^2 l(\theta_0)u + \frac{1}{6}n^{-3/2}\nabla'\left\{ u'\nabla^2 l(\theta^\star)u \right\} u.
$$

where $\theta^\star$ lies in the line segment connecting the points $\theta_0$ and $\theta_0 + n^{-1/2}u$. Standard arguments show that

$$
\begin{aligned}
\mid n^{-1/2}\nabla' l(\theta_0)u \mid &= O_p(1) \| u \|, \\
\frac{1}{2}n^{-1}u'\nabla^2 l(\theta_0)u &= -\frac{1}{2}u'Su + o_p(1), \\
\frac{1}{6}n^{-3/2}\nabla'\left\{ u'\nabla^2 l(\theta^\star)u \right\} u &= o_p(1)
\end{aligned}
$$

where (6) has been used and $S$ has been defined by (7). Similarly, we obtain that

$$
\begin{aligned}
E_2 = -\lambda_n \Bigg( &n^{-1/2} \sum_{j=1}^d \big( q|\beta_j|^{q-1}\mathrm{sign}(\beta_j)u_j I(\beta_j \neq 0) \big) \\
&+ \frac{1}{2}n^{-1} \sum_{j=1}^d q(q-1)\big(|\beta_j|^{q-2}u_j^2 I(\beta_j \neq 0)\big) + o(1) \Bigg)
\end{aligned}
$$

Again, it can be shown that both terms of the above summand are $O(1)$ and $o(1)$ respectively, provided that $\lambda_n/\sqrt{n} \to \lambda_0 \geq 0$ upon recalling (9), (10) and the fact that $\theta_0$ lies in a compact set. For large $\| u \|$ all the terms are small and dominated by $\frac{1}{2}n^{-1}u'\nabla^2 l(\theta_0)u$. Because of Lemma 3.2, $\hat{\theta}^\lambda$ exists with probability 1, hence the theorem follows. If $\lambda_0 = 0$, then the sequence is consistent since $\hat{\theta}^\lambda$ can be chosen independently of $u$.

### Proof of Theorem 3.2

A Taylor expansion shows that

$$0 = \nabla l_p(\hat{\theta}^\lambda) = \nabla l_p(\theta_0) + \nabla^2 l_p(\theta^\star)(\hat{\theta}^\lambda - \theta_0)$$

where $\theta^\star$ lies in the segment between $\theta_0$ and $\hat{\theta}^\lambda$. Rearranging terms in the above expression we obtain

$$
\begin{aligned}
(\hat{\theta}^\lambda - \theta_0) = & -\left\{\nabla^2 l_p(\theta^\star)\right\}^{-1} \nabla l_p(\theta_0) \\
= & -\left\{\nabla^2 l_p(\theta_0)\right\}^{-1} \mathrm{E}\left[\nabla l_p(\theta_0)\right] - \left\{\nabla^2 l_p(\theta^\star)\right\}^{-1} \left\{\nabla l_p(\theta_0) - \mathrm{E}\left[\nabla l_p(\theta_0)\right]\right\} \\
& -\left\{\left\{\nabla^2 l_p(\theta^\star)\right\}^{-1} - \left\{\nabla^2 l_p(\theta_0)\right\}^{-1}\right\} \mathrm{E}\left[\nabla l_p(\theta_0)\right] \\
= & \; D_1 + D_2 + D_3.
\end{aligned}
$$

Then application of the weak law of large numbers together with the central limit theorem shows that

$$
\begin{aligned}
\sqrt{n} D_1 \rightarrow & \quad -\lambda_0 S^{-1} b, \\
\sqrt{n} D_2 \rightarrow & \quad \mathcal{N}_{d+1}(0, V), \\
\sqrt{n} D_3 \rightarrow & \quad 0,
\end{aligned}
$$

as $n \rightarrow \infty$, where the first and third convergence are in probability. These results prove the theorem.

### Proof of Theorem 3.3

Following similar arguments as those of Qin and Zhang (1997), the following representation is obtained by a Taylor expansion

$$
\begin{aligned}
\hat{G}_2^\lambda(t) - G_2(t) \; + & \; \frac{\rho_1}{n}(A_{11}(t), A'_{12}(t)) \left\{-\frac{1}{n}\nabla^2 l(\theta_0) + \frac{\lambda_n}{n}\nabla^2 J(\beta)\right\}^{-1} \mathrm{E}[\nabla l_p(\theta)] \\
= & \; \Lambda(t) - G_2(t) + K(t) + R_n(t),
\end{aligned}
$$

with $\sup_{t \in [-\infty, \infty]} \mid R_n(t) \mid = o_p(n^{-1/2})$,

$$\Lambda(t) = \sum_{i=1}^{2} \sum_{j=1}^{n_i} \frac{I(X_{ij} \le t)}{n_2(1 + \rho_1 \exp(\alpha + \beta' h(x_{ij})))},$$

and

$$
\begin{aligned}
K(t) = & -\frac{\rho_1}{n}(A_{11}(t), A'_{12}(t)) \left\{-\frac{1}{n}\nabla^2 l(\theta_0) + \frac{\lambda_n}{n}\nabla^2 J(\beta)\right\}^{-1} \\
& \times (\nabla l_p(\theta_0) - \mathrm{E}[\nabla l_p(\theta_0)]).
\end{aligned}
$$

Then,

$$\mathrm{E}[\sqrt{n}(\Lambda(t) - G_2(t) + K(t))] = 0$$

and after lengthy calculations

$$
\begin{aligned}
\mathrm{Cov}[\sqrt{n}(\Lambda(t) - G_2(t) + K(t)) \quad ; \quad &\sqrt{n}(\Lambda(s) - G_2(s) + K(s))] = \\
&(1 + \rho_1)(G_2(t \wedge s) - G_2(t)G_2(s)) \\
- \quad &\rho_1(1 + \rho_1)A_{11}(t \wedge s) \\
+ \quad &\rho_1(1 + \rho_1)(A_{11}(t), A'_{12}(t))'A^{-1} \\
\times \quad &(A_{11}(s), A'_{12}(s)).
\end{aligned} \tag{15}
$$

In the above derivation, we use the fact that $\lambda_n/\sqrt{n} \to \lambda_0$ so that

$$
\left\{ -\frac{1}{n}\nabla^2 l(\theta_0) + \frac{\lambda}{n}\nabla^2 J(\beta) \right\}^{-1} - S^{-1} = o_p(1),
$$

where $S$ is recalled by (7). Therefore the central limit theorem and the Cramer–Wold device imply that the finite dimensional distributions of the process $\sqrt{n}(\Lambda(t) - G_2(t) + K(t))$ converge weakly to those of a Gaussian process with mean 0 and covariance function given by (15). In order to prove the conclusion of the theorem, it is sufficient to show $\{\sqrt{n}(\Lambda(t) - G_2(t) + K(t)), t \in [-\infty, \infty]\}$ is tight in $D[-\infty, \infty]$–a fact which can be shown by using tightness criteria.